\documentclass[reqno]{amsart}

\usepackage{amsfonts}
\usepackage{amssymb}
\usepackage{times}
\usepackage{mathptmx}

\usepackage{amsmath}

\newcommand{\diag}{\mathop{\rm diag}\nolimits}
\newtheorem{theorem}{Theorem}[section]
\newtheorem{lemma}{Lemma}[section]

\begin{document}



\title{Inequalities with Determinants of Perturbed Positive Matrices}


\author{Ivan Matic}

\address{Department of Mathematics, Baruch College, CUNY,
One Bernard Baruch Way, New York, NY 10010, USA}

\maketitle

\begin{abstract}
 We prove two inequalities regarding the ratio $\det(A+D)/\det A$ of the determinant of a positive-definite matrix $A$ and the determinant of its perturbation $A+D$. In the first problem, we study the perturbations that happen when positive matrices are added to diagonal blocks of the original matrix. In the second problem, the perturbations are added to the inverses of the matrices.



\end{abstract}


\section{Introduction}
\label{in}

Given $k$ complex square matrices $B_1$, $\dots$, $B_k$ of format $n_1\times n_1$, $\dots$, $n_2\times n_2$, $\dots$, $n_k\times n_k$, let us denote by $\diag(B_1, \dots, B_k)$ the matrix of the format $(n_1+\cdots+n_k)\times (n_1+\cdots+n_k)$ whose main diagonal blocks are $B_1$, $\dots$, $B_k$ and all other entries are $0$. In other words:
\begin{eqnarray}\nonumber\diag(B_1,\dots, B_k)=\left[\begin{array}{cccc}B_1&0&\cdots&0\\0&B_2&\cdots&0\\ 0&0&\ddots&0\\ 
0&0&\cdots&B_k\end{array}\right].\end{eqnarray} 
 Given two vectors $u,v\in\mathbb C^k$ such that $u=\langle u_1,\dots, u_k\rangle$ and $v=\langle v_1, \dots, v_k\rangle$, we define their inner product $\langle u,v\rangle= \sum_{i=1}^k u_i\overline{v_i}$. For a complex $n\times m$ matrix $R$, we use $R^*$ to denote its adjoint matrix. In other words, $R^*$ is the transpose of the complex conjugate of $R$, and for $u\in\mathbb C^n$ and $v\in\mathbb C^m$ the following is satisfied: $\langle Ru,v\rangle =\langle u,R^*v\rangle$. A square matrix $A$ is self-adjoint if $A^*=A$.

The self-adjoint matrix $A$ of format $n\times n$ is called positive  (or positive definite) if $\langle Ax,x\rangle >0$ for each non-zero vector $x\in\mathbb C^n$. If the strict inequality is replaced by $\geq$, the matrix is called non-negative (or positive semi-definite).  If $A$ and $B$ are two square matrices of the same format, we will write $A\geq B$ (resp. $A>B$) if $A-B\geq 0$ (resp. $A-B>0$).

For $n\in\mathbb N$ we will denote by $I_n$ the $n\times n$ identity matrix.  The subscript $n$ will be omitted when there is no danger of ambiguity.

We will prove the following two inequalities regarding positive matrices with complex entries.

\begin{theorem}\label{sub_02} Assume that $k\in\mathbb N$ and that $n_1$, $\dots$, $n_k$ are positive integers. Assume that $(C_i)_{i=1}^k$ and $(D_i)_{i=1}^k$ are two sequences of positive matrices such that for each $i\in\{1,2,\dots, k\}$ the matrices $C_i$ and $D_i$ are of format $n_i\times n_i$. Assume that $C$ is a positive matrix whose diagonal blocks are $C_1$, $\dots$, $C_k$. The following inequality holds:
\begin{eqnarray}\label{main_inequality}\frac{\det\left(C+\diag(D_1,\dots, D_k)\right)}{\det C}\geq \frac{\det(C_1+D_1)}{\det C_1}\cdots\frac{\det(C_k+D_k)}{\det C_k}.
\end{eqnarray}
\end{theorem}

\begin{theorem}\label{sub_01}  Assume that $k\in\mathbb N$ and that $n_1$, $\dots$, $n_k$ are positive integers. Assume that $(C_i)_{i=1}^k$ and $(D_i)_{i=1}^k$ are two sequences of positive matrices such that for each $i\in\{1,2,\dots, k\}$ the matrices $C_i$ and $D_i$ are of format $n_i\times n_i$. Assume that $C$ is a positive matrix such that the diagonal blocks of $C^{-1}$ are  $C_1^{-1}$, $\dots$, $C_k^{-1}$. The following inequality holds:
\begin{eqnarray}\label{second_inequality}\frac{\det\left(C+\diag(D_1,\dots, D_k)\right)}{\det C}\leq \frac{\det(C_1+D_1)}{\det C_1}\cdots\frac{\det(C_k+D_k)}{\det C_k}.
\end{eqnarray}
\end{theorem}

The two inequalities presented in this paper have the flavor of Fischer's determinantal inequality, although in (\ref{main_inequality}) the sign is reversed. 
An inequality related to our results, which features quotients of perturbed matrices, has been established previously  \cite{lewent_mlin}. For  refinements of Fischer-type inequalities with singular values, the reader is referred to \cite{drury_mlin} and \cite{fu_he}. 
After taking the logarithms of left and right sides of the inequality  (\ref{main_inequality}), one obtains $$\varphi (C,\diag(D_1,D_2))\leq \varphi(C_1,D_1)+\varphi(C_2,D_2),$$ where $\varphi(X,Y)=\log\det(X+Y)-\log\det(X)$. Similar inequalities are known to hold for concave  functions $\varphi$, and such results can be found in  \cite{aujla_bourin}.

The proof of the first theorem relies on  Lemma \ref{decreasing_fraction} which is established using the Grothendieck's determinantal inequality. The lemma implies that the map $U\mapsto \det(U+D)/\det(U)$ is operator-decreasing. Several results about operator-monotone functions are available in \cite{hansen_frank}. Generalizations and improvements of the Grothendieck's inequality have been established in \cite{bourin_lee_mlin} and \cite{seiler_simon} and they have been used in the past to prove results regarding block matrices.

The inequality (\ref{second_inequality}) can be used to establish super-additivity for functions of diffusions in random environments. More precisely, let us consider one-dimensional Brownian motion $Z$, and let $W$ be another Brownian motion independent on $Z$. Define \begin{eqnarray} \label{ld_ineq} f(t)=\log\mathbb E\left[\exp\left(-\int_0^t\left|W(Z(s))\right|^2\,ds\right)\right],\end{eqnarray} where $\mathbb E$ denotes the expected value with respect to the Brownian motion  $W$. We will now illustrate that $f(t_1+t_2)\geq f(t_1)+f(t_2)$ is a special case of the inequality (\ref{second_inequality}). Assume that $0=s_0<s_1<\cdots <s_n=t_1$ is the partition of the interval $[0,t_1]$ into $n$ sub-intervals of length $\lambda_1$. Similarly, let $t_1=s_{n}<s_{n+1}<\cdots<s_{n+m}=t_1+t_2$ be the partition of the interval $[t_1,t_1+t_2]$ into $m$ intervals of length $\lambda_2$. Let us denote \begin{eqnarray*}W_i&=&W(Z(s_i)), \\ \overrightarrow{w_1}&=&\langle W_1,\dots, W_n\rangle,\\ \overrightarrow{w_2}&=&\langle W_{n+1},\dots, W_{n+m}\rangle, \mbox{ and }\\\overrightarrow w&=&\langle W_1,\dots, W_{n+m}\rangle.\end{eqnarray*} Then $\sum_{i=1}^n|W_i|^2\lambda_1=\langle \lambda_1I\overrightarrow{w_1},\overrightarrow{w_1}\rangle$, $\sum_{i=n+1}^{n+m}|W_i|^2\lambda_2=\langle \lambda_2I\overrightarrow{w_2},\overrightarrow{w_2}\rangle$, and $$\sum_{i=1}^n|W_i|^2\lambda_1+\sum_{i=n+1}^{n+m}|W_i|^2\lambda_2=\left\langle \diag(\lambda_1I_n,\lambda_2I_m)\overrightarrow w,\overrightarrow w\right\rangle.$$
If we fix the Brownian motion $Z$, then  $[W(Z(s_1)),\dots, W(Z(s_2))]$ is a multivariate Gaussian random variable and as such it has a covariance matrix $C$. Denote by $C_1$ and $C_2$ the covariance matrices of $[W(Z(s_1)),\dots, W(Z(s_n))]$ and $[W(Z(s_{n+1})), 
\dots, W(Z(s_{n+m}))]$. Then $C_1$ and $C_2$ are the diagonal blocks of $C$. Moreover,
\begin{eqnarray*}\mathbb E\left[\exp(-\langle \lambda_1I\overrightarrow{w_1},\overrightarrow{w_1}\rangle)\right]&=&\frac1{M\sqrt{\det C_1}}\int e^{-\langle \lambda_1\overrightarrow{w_1},\overrightarrow{w_1}\rangle-\frac12\langle C_1^{-1}\overrightarrow{w_1},\overrightarrow{w_1}\rangle}\,d\overrightarrow{w_1}\\
&=&\frac1{M\sqrt{\det C_1}}\int e^{-\left\langle \left(C_1^{-1}+\lambda_1I\right)\overrightarrow{w_1},\overrightarrow{w_1}\right\rangle}\,d\overrightarrow{w_1}\\
&=&\frac{\sqrt{ \det C_1^{-1}}}{\sqrt{\det\left(C_1^{-1}+\lambda_1I\right)}},
\end{eqnarray*}
where $M$ is a normalizing constant. 
We obtain analogous equalities for the quantities $\mathbb E\left[\exp(-\langle \lambda_2I\overrightarrow{w_2},\overrightarrow{w_2}\rangle)\right]$ and 
$\mathbb E\left[\exp(-\langle \diag(\lambda_1I_n,\lambda_2I_m)\overrightarrow{w},\overrightarrow{w}\rangle)\right]$. 
The inequality (\ref{second_inequality}) implies that \begin{eqnarray*}&&\log\mathbb E\left[\exp\left(-\left\langle\diag(\lambda_1I_n,\lambda_2I_m)\overrightarrow w,\overrightarrow w\right\rangle\right)\right]\\&\geq&
\log\mathbb E\left[\exp\left(-\langle \lambda_1I\overrightarrow{w_1},\overrightarrow{w_1}\rangle\right)\right]+\log\mathbb E\left[\exp\left(-\langle \lambda_2I\overrightarrow{w_2},\overrightarrow{w_2}\rangle\right)\right].\end{eqnarray*}
Taking the limit as $\lambda_1,\lambda_2\to 0$ we obtain $f(t_1+t_2)\geq f(t_1)+f(t_2)$.

This technique is potentially useful for establishing large deviations  for  random processes with drifts. Sub-additive properties are known to hold for  killed Brownian motions in random environments \cite{sznitmanbook, zerner}.  However, in the case of drifts introduced to random diffusions, no analogous results have yet been established. If a drift is assumed to be a multivariate Gaussian process, a possible approach is to express the large deviation probabilities in terms of determinants.   However, there is still work to be done to transform the general case of sub-additive inequalities into the language  of their covariance matrices \cite{d_proc_re}.

\section{Theorems from literature that are used in the proofs}

We will start with listing the known theorems that we will use to establish the inequalities. For the derivations of the results presented in this section, the reader is referred  to  \cite{horn_johnson}.
\begin{theorem} (Weyl's inequality) \label{determinants}
If $A$ and $B$ are non-negative matrices such that $A\geq B$ then $\det A\geq \det B$. If both of them are invertible then $A^{-1}\leq B^{-1}$.
\end{theorem}

\begin{theorem} (Block inverse theorem) \label{block_inverse} 
Let $A$ and $D$ be square $n\times n$ and $m\times m$ matrices respectively. Assume that $B$ and $C$ are matrices of the formats $n\times m$ and $m\times n$ and assume that $A$ and 
$S_A=D-CA^{-1}B$ are invertible. Then 
$$\left[\begin{array}{cc}A&B\\C&D\end{array}\right]^{-1}=
\left[\begin{array}{cc}A^{-1}+A^{-1}BS_A^{-1}CA^{-1}&-A^{-1}BS_A^{-1}\\
-S_A^{-1}CA^{-1}&S_A^{-1}\end{array}\right].
$$ 
Similarly, if $D$ and $S_D=A-BD^{-1}C$ are invertible, then
$$\left[\begin{array}{cc}A&B\\C&D\end{array}\right]^{-1}= \left[\begin{array}{cc}S_D^{-1}&-
S_D^{-1}BD^{-1}
\\
-D^{-1}CS_D^{-1}
&
D^{-1}+D^{-1}CS_D^{-1}BD^{-1}
\end{array}\right].
$$
\end{theorem}

The matrices $S_A$ and $S_D$ are called Schur complements of $A$ and $D$.

 The following two consequences of the previous result are known as the Woodbury's matrix identity and the Fischer's inequality.
\begin{theorem}\label{sasd}
Let $A$, $B$, $C$, $D$, $S_A$, $S_D$ be as in Theorem \ref{block_inverse}. Then
\begin{eqnarray*}(A-BD^{-1}C)^{-1}&=&A^{-1}+A^{-1}BS_A^{-1}CA^{-1},
\end{eqnarray*}
provided that the inverses are defined. Moreover if the matrix $M=\left[\begin{array}{cc}A&B\\C&D\end{array}\right]$ is positive, then $S_A$ and $S_D$ are positive. 
\end{theorem}
\begin{theorem}\label{fischer_ineq} 
Let $A$, $B$, $C$, $D$, $S_A$, $S_D$ be as in Theorem \ref{block_inverse}. Let $M=\left[\begin{array}{cc}A&B\\C&D\end{array}\right]$. Then
\begin{eqnarray*}\det M=\det A\cdot \det S_A.
\end{eqnarray*}
If the matrix $M$ is positive then $\det S_A\leq \det D$, and $\det M\leq \det A\cdot \det D$.
\end{theorem}

\begin{theorem} (Grothendieck \cite{grothendieck_55}) If $A$ and $B$ are non-negative symmetric matrices of the format $n\times n$ and $I$ the $n\times n$ identity matrix then 
\begin{eqnarray}\label{gss_ineq}\det(I+A+B)\leq \det(I+A)\det(I+B).\end{eqnarray}
\end{theorem}

\section{Proofs of Theorems \ref{sub_02} and \ref{sub_01}}
We will start by proving Theorem \ref{sub_01} since it is easier to prove than Theorem \ref{sub_02}.

\noindent{\bf Proof of Theorem \ref{sub_01}.} 
Let us denote $B_i=C_i^{-1}$ for $i\in\{1,2,\dots, k\}$ and $B=C^{-1}$. Using the multiplicative property of determinants we transform the inequality (\ref{second_inequality}) into equivalent one: 
\begin{eqnarray}\label{equivalent_form}
\det(I+B\diag(D_1,\dots, D_k))\leq \det(I+B_1D_1)\cdots \det (I+B_kD_k).
\end{eqnarray}
Since the matrices $D_1$, $\dots$, $D_k$ are positive we have that each of them has a square root. In other words, for each $i$, there exists a unique positive matrix $\sqrt{D_i}$ that commutes with $D_i$ and satisfies $D_i=\sqrt{D_i}\cdot \sqrt{D_i}$. Clearly, the matrix $\diag\left(\sqrt{D_1},\dots, \sqrt{D_k}\right)$ is the square root of $\diag(D_1,\dots, D_k)$. 

Applying the Sylvester's determinant identity $\det(I+XY)=\det(I+YX)$ to the matrices $X=B\sqrt{\diag(D_1,\dots, D_k)}$ and $Y=\sqrt{\diag(D_1,\dots,D_k)}$ we transform the left-hand side of (\ref{equivalent_form}) into:
\begin{eqnarray*} &&\det\left(I+B\diag(D_1,\dots, D_k)\right)\\&=&\det\left(I+\diag\left(\sqrt{D_1},\cdots, \sqrt{D_k}\right)\cdot B\cdot \diag\left(\sqrt{D_1},\cdots, \sqrt{D_k}\right)\right).
\end{eqnarray*}
Similarly, the right-hand side of (\ref{equivalent_form}) is:
\begin{eqnarray*} &&
\det(I+B_1D_1)\cdots \det(I+B_kD_k)\\&=&\det\left(I+\sqrt{D_1}B_1\sqrt{D_1}\right)\cdots\left(I+\sqrt{D_k}B_k\sqrt{D_k}\right).
\end{eqnarray*}
We will use the induction on $k$ to prove the inequality (\ref{equivalent_form}). Let us start with $k=2$ and  assume that $B=\left[\begin{array}{cc}B_1&R\\R^*&B_2\end{array}\right]$ for an $n_1\times n_2$ matrix $R$. Elementary calculations imply:
\begin{eqnarray*}&&\det\left(I+\diag\left(\sqrt{D_1},\sqrt{D_2}\right)\cdot B\cdot\diag\left(\sqrt{D_1}, \sqrt{D_2}\right)\right)\\&=&
\det\left[\begin{array}{cc}I+\sqrt{D_1}B_1\sqrt{D_1} & \sqrt{D_1}R\sqrt{D_2}\\ \sqrt{D_1}R^*\sqrt{D_1} & I+\sqrt{D_2}B_2\sqrt{D_2}\end{array}\right].
\end{eqnarray*}
We can now use Theorem \ref{fischer_ineq} to conclude 
that \begin{eqnarray*}&&\det\left(I+\diag\left(\sqrt{D_1},\sqrt{D_2}\right)\cdot B\cdot\diag\left(\sqrt{D_1}, \sqrt{D_2}\right)\right)\\&\leq& \det\left(I+\sqrt{D_1}B_1\sqrt{D_1}\right)\cdot \det\left(I+\sqrt{D_2}B_2\sqrt{D_2}\right)\\
&=& \det\left(I+D_1B_1\right)\cdot\det\left(I+D_2B_2\right).\end{eqnarray*}
Therefore the inequality (\ref{equivalent_form}) is established for $k=2$.

Assume now that $k\geq 3$ and that the inequality (\ref{equivalent_form}) is true for $k-1$. Assume that $D_1$, $\dots$, $D_k$, and $B_1$, $\dots$, $B_k$ are positive matrices of formats $n_1\times n_1$, $\dots$, $n_k\times n_k$. Assume that $B_2'$ is the sub-matrix of the matrix $B$ obtained by removing the first $n_1$ rows and first $n_1$ columns. 
According to the induction hypothesis we have 
\begin{eqnarray}\label{eqn_km1}
\det\left(I+B'_2\diag\left(D_2,\dots, D_k\right)\right)&\leq& \det(I+B_2D_2)\cdots\det(I+B_kD_k).
\end{eqnarray}
Let us denote $D'_2=\diag(D_2,\dots, D_k)$. We can write $\diag(D_1,\dots, D_k)=\diag(D_1,D_2')$. Applying the inequality (\ref{equivalent_form}) with $k=2$ we obtain
\begin{eqnarray}\label{eqn_km1_2}
\det\left(I+B\diag(D_1,D_2')\right)&\leq&\det (I+B_1D_1)\cdot \det(I+B_2'D_2').
\end{eqnarray}
The inequalities (\ref{eqn_km1}) and (\ref{eqn_km1_2}) together imply the inequality (\ref{equivalent_form}). This completes the proof of Theorem \ref{sub_01}. \hfill $\Box$

\vspace{0.2cm}

In order to prove Theorem \ref{sub_02} we will need the following lemma. 

\begin{lemma}\label{decreasing_fraction} Assume that $U\geq V$ and $D$ are $n\times n$ non-negative matrices such that $U$ and $V$ are invertible. Then the following inequality holds: $$\frac{\det (V+D)}{\det V}\geq \frac{\det(U+D)}{\det U}.$$
\end{lemma}

\noindent{\bf Proof.}
 The matrix $V^{-1}$ is positive and as such it has a positive square root. Let us denote it by $V^{-\frac12}$. Assume that $U=V+W$ for some non-negative matrix $W$. The required inequality is equivalent to
\begin{eqnarray*}\det(V+W)\det(V+D)&\geq& \det V\cdot \det(V+W+D).\end{eqnarray*}
We now multiply both left and right side of the previous inequality by $\left[\det\left(V^{-\frac12}\right)\right]^4$:
\begin{eqnarray*}
&& 
\det\left(V^{-\frac12}\right)\det(V+W)\det\left(V^{-\frac12}\right)\cdot \det\left(V^{-\frac12}\right)\det(V+D)\det\left(V^{-\frac12}\right)\\&&\geq
\det\left(V^{-\frac12}\right)\det(V+W+D)\det\left(V^{-\frac12}\right). \end{eqnarray*}
The last inequality is equivalent to:
\begin{eqnarray*}
&&
\det\left(I+V^{-\frac12}WV^{-\frac12}\right)\cdot \det\left(I+V^{-\frac12}DV^{-\frac12}\right)\\&\geq & 
\det\left(I+V^{-\frac12}WV^{-\frac12}+V^{-\frac12}DV^{-\frac12}\right).\end{eqnarray*}
The last inequality can be derived by applying (\ref{gss_ineq}) to the positive definite matrices $A=V^{-\frac12}WV^{-\frac12}$ and $B=V^{-\frac12}DV^{-\frac12}$.
\hfill $\Box$

\vspace{0.2cm}

\noindent{\bf Proof of Theorem \ref{sub_02}.} We will first prove the theorem for the case $k=2$. Assume that $C=\left[\begin{array}{cc}C_1&R\\R^*&C_2\end{array}\right]$ for some matrix $R$ of the format $n_1\times n_2$. From Theorem \ref{fischer_ineq} we conclude that the required inequality is  equivalent to
\begin{eqnarray*}&&\frac{\det (C_1+D_1)\det\left[(C_2+D_2)-R^*(C_1+D_1)^{-1}R\right]}{\det C_1\det\left(C_2-R^*C_1^{-1}R\right)}\\&\geq &
\frac{\det(C_1+D_1)\cdot \det(C_2+D_2)}{\det C_1\cdot \det C_2}.
\end{eqnarray*}
This inequality can be re-written as $$\frac{\det\left[(C_2+D_2)-R^*(C_1+D_1)^{-1}R\right]}{\det(C_2+D_2)}\geq\frac{\det\left(C_2-R^*C_1^{-1}R\right)}{\det C_2}.
$$
Let us denote $V=C_2-R^*(C_1+D_1)^{-1}R$. We can prove that $V>0$ by applying Theorem \ref{sasd} to the matrix $\tilde C=\left[\begin{array}{cc}C_1+D_1&R\\R^*&C_2\end{array}\right]$. The theorem requires the positivity of $\tilde C$, and this is satisfied since $\tilde C=C+\diag(D_1,0)$. Let $U=C_2$. Moreover, $U-V=R^*(C_1+D_1)^{-1}R\geq 0$, therefore we can apply Lemma \ref{decreasing_fraction} to matrices $U$ and $V$ to obtain:
$$\frac{\det\left[(C_2+D_2)-R^*(C_1+D_1)^{-1}R\right]}{\det(C_2+D_2)}\geq\frac{\det\left[C_2-R^*(C_1+D_1)^{-1}R\right]}{\det C_2}.
$$
From $C_1+D_1\geq C_1$ we have $(C_1+D_1)^{-1}\leq C_1^{-1}$. Therefore \begin{eqnarray*}R^*(C_1+D_1)^{-1}R&\leq& R^*C_1^{-1}R,\quad\mbox{ and }\\
C_2-R^*(C_1+D_1)^{-1}R&\geq&  C_2-R^*C_1^{-1}R.\end{eqnarray*} Theorem \ref{determinants} now implies   $$\det\left(C_2-R^*(C_1+D_1)^{-1}R\right)\geq
\det \left(C_2-R^*C_1^{-1}R\right)$$ which completes the proof of Theorem \ref{sub_02} when $k=2$.

We will use induction to finish the proof for general $k\in\mathbb N$. Assume that $k\geq 3$ and that the statement is true for $k-1$. We will now prove the inequality for  matrices $C_1$, $\dots$, $C_k$, $D_1$, $\dots$, $D_k$.
Let us denote by $C_2'$ the sub-matrix of the matrix $C$ obtained by removing its first $n_1$ rows and first $n_1$ columns. The matrix $C$ can be regarded as a block matrix with diagonal blocks $C_1$ and $C_2'$. Similarly, for $D_2'=\diag(D_2,\dots, D_k)$ we have that $\diag(D_1,\dots, D_k)=\diag (D_1,D_2')$. Using the induction hypothesis we obtain
\begin{eqnarray}\label{ineq_01}\frac{\det\left(C_2'+\diag(D_2,\dots, D_k)\right)}{\det C_2'}&\geq &\frac{\det(C_2+D_2)}{\det C_2}\cdots\frac{\det(C_k+D_k)}{\det C_k}.\end{eqnarray}
Using the inequality established for $k=2$ we conclude 
\begin{eqnarray}\nonumber \frac{\det\left(C+\diag(D_1,\dots, D_k)\right)}{\det C}&=&\frac{\det\left(C+\diag(D_1, D'_2)\right)}{\det C}\\ \label{ineq_02}&\geq &
\frac{\det(C_1+D_1)}{\det C_1}\cdot\frac{\det(C_2'+D_2')}{\det C_2'}.\end{eqnarray}
The inequalities (\ref{ineq_01}) and (\ref{ineq_02}) imply the desired result. \hfill$\Box$

\vspace{0.2cm}

\section{Effects of perturbations by arbitrary positive matrices}

In this section we will show that neither of the above inequalities can be generalized to allow for $\diag(D_1,\dots, D_k)$ to be replaced by an arbitrary positive matrix whose diagonal blocks are $D_1$, $\dots$, $D_k$. We first present an example that illustrates the case in which the reverse inequality occurs in (\ref{main_inequality}) under previously mentioned generalization. Take $$C=\left[\begin{array}{cc}10& 2\\2& 5\end{array}\right] \quad\mbox{and}\quad D=\left[\begin{array}{cc}2& 1\\ 1& 1\end{array}\right].$$ Then $C_1=10$, $C_2=5$, $D_1=2$, $D_2=1$ and a simple calculation shows that the left-hand side of (\ref{main_inequality}) corresponds to $$\frac{\det (C+D)}{\det C}=\frac{63}{46},$$ while the right-hand side is $$\frac{\det(C_1+D_1)}{\det C_1}\cdot \frac{\det(C_2+D_2)}{\det C_2}=\frac{72}{50}=\frac{36}{25}>\frac{63}{46}.$$

To see a counter-example to (\ref{second_inequality}) if we allow for $\diag(D_1,\dots, D_k)$ to be replaced by general $D$, we consider 
$$C=\left[\begin{array}{cc}2&-2\\-2&4\end{array}\right] \quad\mbox{and}\quad D=\left[\begin{array}{cc}1&1\\1&2\end{array}\right].
$$

We now have $D_1=1$ and $D_2=2$. In order to determine $C_1$ and $C_2$ we first find $$C^{-1}=\left[\begin{array}{cc}1&\frac12\\ \frac12&\frac12\end{array}\right]$$
which implies that $C_1=1$ and $C_2=2$. It is now easy to find that the left-hand side of (\ref{second_inequality}) corresponds to $$\frac{\det(C+D)}{\det C}=\frac{17}4>4=
\frac{\det(C_1+D_1)}{\det C_1}\cdot \frac{\det(C_2+D_2)}{\det C_2}.$$ Thus, the inequality (\ref{second_inequality}) does not always hold if the matrix $\diag(D_1,\dots, D_k)$ is replaced with a  positive matrix $D$ of a more general form.

\section{Acknowledgments} The author expresses gratitude to Fraydoun Rezakhanlou for discussions and ideas regarding the problems in this paper. The author would also like to thank the anonymous referee for helpful suggestions and comments.




\begin{thebibliography}{00}



\bibitem{aujla_bourin} J.S. Aujla, J-C Bourin, {\em Eigenvalue inequalities for convex and log-convex functions}, Linear Algebra Appl. 424 (2007), no. 1, 25--35.

\bibitem{bourin_lee_mlin} J-C Bourin, E-Y Lee, M. Lin, {\em Positive matrices partitioned into a small number of Hermitian blocks}, Linear Algebra Appl. 438 (2013), no. 5, 2591--2598.

\bibitem{drury_mlin} S. Drury, M. Lin, {\em Reversed Fischer determinantal inequalities}, to appear in Linear Multilinear Algebra (2014).

\bibitem{fu_he} X. Fu, C.He, {\em On some Fischer-type determinantal inequalities for accretive-dissipative matrices}, J. Inequal. Appl. 2013, 2013:316, 7 pp. 

\bibitem{grothendieck_55} A. Grothendieck,  {\em R\'arrangements de functions et in\'egalit\'es de convexit\'e dans le alg\`ebres de von Neumann d'une trace} (mimeographed notes), in 
``S\'eminaire Bourbaki'' (1955)  pp. 113.01--113.13

\bibitem{hansen_frank} F. Hansen, {\em WYD-like skew information measures}, J. Stat. Phys. 151 (2013), no. 5, 974--979. 
\bibitem{lewent_mlin} M. Lin, {\em A Lewent type determinantal inequality}, Taiwanese J. Math. 17 (2013), no. 4, 1303--1309. 

\bibitem{horn_johnson} R.A. Horn, C.R. Johnson, {\em Matrix analysis}, Cambridge University Press, 1990.

\bibitem{d_proc_re} I. Matic, {\em Large deviations for processes in random environments with jumps}, Electron. J. Probab. 16 (2011), no. 87, pp. 2406--2438.

\bibitem{seiler_simon} E. Seiler, B. Simon,  {\em An inequality among determinants}, Proc. Nat. Acad. Sci. USA, Vol 72 (1975),  No. 9, pp. 3277--3278.
\bibitem{sznitmanbook} A-S Sznitman,  {\em Brownian Motion, Obstacles and Random Media}, Springer Monographs in Mathematics, 1998.



\bibitem{zerner} M.P.W. Zerner,  {\em Directional decay of the Green's function for a random nonnegative potential on $\mathbb Z^d$}, Ann. Appl. Probab., 8 (1997),  246--280.


\end{thebibliography}


\end{document}